\documentclass[a4paper,10pt]{amsart}
\usepackage{amsfonts}
\usepackage{amsthm}
\usepackage{amssymb}
\usepackage{amstext}
\usepackage{amsmath}
\usepackage{wasysym}
\usepackage{tikz}

\usepackage{graphics}
\usepackage{epsfig}
\usepackage{amscd}
\usepackage[all]{xy}
\usepackage{latexsym}
\usepackage{graphicx}
\usepackage{multirow}
\usepackage{geometry}

\usepackage{fancyhdr}
\usepackage{lipsum}
\usepackage{indentfirst}
\usepackage{graphicx}
\usepackage{textcomp}
\usepackage{faktor}
\usepackage{newlfont}
\usepackage{latexsym}
\usepackage{mathtools}

\newcommand{\Q}{\mathbb{Q}}
\newcommand{\cE}{\mathcal{E}}

\newcommand{\cor}[1]{\left\langle #1 \right\rangle }

\theoremstyle{plain}                    
\newtheorem{teo}{Theorem}[section]      
\theoremstyle{definition}               
\newtheorem{ese}{Example}      
\theoremstyle{remark}                   
     
\def\SS{\mathcal{S}}
\def\EE{\mathcal{E}}

\def\J{\mathcal{J}}

\numberwithin{equation}{section}

\begin{document}

\title[Harer-Zagier formula via Fock space]
      {Harer-Zagier formula via Fock space}
      
\author[D.~Lewa\'{n}ski]{D.~Lewa\'{n}ski}
\address{D.~L.: Max Planck Institut für Mathematik, Vivatsgasse 7, 53111 Bonn, Germany.}
\email{ilgrillodani@mpim-bonn.mpg.de}


\maketitle

Let $\epsilon_g(d')$ be the number of ways of obtaining a genus $g$ Riemann surface by identifying in pairs the sides of a $(2d')$-gon.
The goal of this note is to give a short proof of the following theorem.
\begin{teo}[Harer-Zagier formula, \cite{HZ}, 1986]\label{thm:hz}
\begin{equation}
  \epsilon_g(d')= 
 \frac{(2d'-1)!! 2^{d'-2g}}{(d' - 2g + 1)!} [u^{2g}] \left[ \frac{u}{\sinh(u)}\right]^{2} \left[ \frac{u }{\tanh(u)} \right]^{d'}.
\end{equation}
\end{teo}
\noindent
This formula was needed in the same paper \cite{HZ} as key combinatorial fact to compute the celebrated formula for the Euler characteristic of the moduli space of curves of genus $g$. Nine years after the original paper, a proof of this formula via representation theory was given in \cite{DON2}. Later on, it was proved again using different techniques (see, e.g., \cite{AS, Cha, CD, Pittel}). Our proof uses operators in the semi-infinite wedge formalism. The operator we need is available in the literature, but corresponds to a different enumerative problem. However, this enumerative problem is known to be related to $\epsilon_g(d')$ via a certain correspondence in the group algebra of the symmetric group due to Jucys.

\subsection{Hurwitz numbers as vacuum expectations} The pioneering paper of Okounkov in 2000 \cite{Ok00} expresses simple Hurwitz numbers as vacuum expectation via the so called semi-infinite wedge formalism (for an introduction on the topic we recommend \cite{Joh10}). Since then, many more (generating series of) different kinds of Hurwitz numbers have been related to vacuum expectations of the type
$$
\cor{\mathcal{E}_{a_1}(z_1) \mathcal{E}_{a_2}(z_2) \cdots \mathcal{E}_{a_n}(z_n) }, \qquad\qquad  \qquad \sum_{i=1}^n a_i = 0,
$$
where $a_i$ are integer numbers called \textit{energies}, and $z_i$ are complex variables. Such expressions vanish whenever $a_1$ is negative or $a_n$ is positive \textemdash\, the general strategy to evaluate them is to commute the $\mathcal{E}$ operators exploiting this vanishing property, until a single energy zero operator $\mathcal{E}_0(z)$ is left:
\begin{equation}\label{eq:commE}
  \left[\cE_a(z),\cE_b(w)\right] = 
  \begin{cases}
  a \delta_{a+b, 0}  &\text{ if } z = w = 0, \\
  \varsigma( aw - bz) \mathcal{E}_{a+b}(z+w) &\text{ otherwise,}
  \end{cases}
  \qquad \qquad \qquad
  \cor{\mathcal{E}_a(z)} = 
  \frac{\delta_{a,0}}{\varsigma(z)},
\end{equation}
with $\varsigma(z) := 2\sinh(z/2) = e^{z/2} - e^{-z/2}$. For the purpose of this note we will actually only need the last computational rule. We will also use the function $\mathcal{S}(z) = \varsigma(z)/z$.

\subsection{Recalling four facts} Our proof consists of an elementary computation, whose first equation follows from four known facts.
\begin{itemize}
\item[\textit{Fact 1.}]  Standard topological arguments (see, e.g., \cite{DM}) show that
$$ \frac{\epsilon_g(d')}{(2d')} = h_{g; (2d')}^{Gr., 2}$$
where $h_{g; (2d')}^{\text{Gr.}, 2}$ is the weighted enumeration of ramified coverings of the Riemann sphere $\mathbb{C}\mathbb{P}^1$ of degree $2d'$ by a genus $g \geq 0$ surface only ramified over $0, 1 , \infty \,\in \,\mathbb{C}\mathbb{P}^1$ (also known as Grothendieck dessins d'enfant), whose ramification profiles over $0$ and $\infty$ are given respectively by the cycle types $(2d')$ and $(2, 2, \dots, 2)$, whereas the cycle type $\mu \vdash 2d'$ over $1$ is arbitrary (also known as one-part $2$-hypermaps, or ribbon graphs with a single face). 
\item[\textit{Fact 2.}] Jucys Correspondence \cite{J} implies that
$$ h_{g; (2d')}^{Gr., 2} = h_{g; (2d')}^{<, 2},$$
where $h_{g; (2d')}^{<, 2}$ is the weighted enumeration of ramified coverings of the Riemann sphere $\mathbb{C}\mathbb{P}^1$ of degree is $2d'$ by a genus $g \geq 0$ surface ramified over $0, \infty \,\in \,\mathbb{C}\mathbb{P}^1$, with ramification profiles given respectively by the cycle types $(2d')$ and $(2, 2, \dots, 2)$, and $m$ further simple ramifications $(x_i \; y_i)_{i = 1, \dots, m} \, \in \, \mathfrak{S}_{2d'}$ written such that $x_i<y_i$ and subject to the extra \textit{strictly monotone} condition $y_i < y_{i+1}$. This can be found, e.g. in \cite{ALS, HO}. We prove it once more in Appendix \ref{app} for reader's convenience.
\item[\textit{Fact 3.}] 
The vacuum expectation operator for strictly monotone ramifications $\mathcal{D}^{(\sigma)}(u)$ was derived in \cite{ALS}, Proposition 5.2, and it implies that
\begin{equation*}
\sum_{g=0} h_{g; (2d')}^{<, 2} u^{2g - 1 + d'}= \frac{1}{(2d')} \cor{e^{\cE_2(0)/2} \mathcal{D}^{(\sigma)}(u) \cE_{-2d'}(0) \left( \mathcal{D}^{(\sigma)}(u) \right)^{-1} e^{-\cE_2(0)/2}},
\end{equation*}
where
$
\mathcal{D}^{(\sigma)}(u) = \exp\left(\frac{\cE_0(-u^2\partial_u))}{\varsigma(u^2\partial_u)} + [w^1]\cE_0(w) \right)\! \log(u).
$
\item[\textit{Fact 4.}]
The nested conjugations above were computed in a more general setting in \cite{KLS}, Lemma 4.4. The computation is an elementary application of the Campbell-Baker-Hausdorff formula $e^A B e^{-A} = \sum_{t=0}\frac{1}{t!} [A, \dots [A,[A,B]]\dots ]$, since $\mathcal{D}^{(\sigma)}(u)$ itself is an exponential operator. It implies that
\begin{multline*}
e^{\cE_2(0)/2} \mathcal{D}^{(\sigma)}(u) \cE_{-2d'}(0) \left( \mathcal{D}^{(\sigma)}(u) \right)^{-1}\!\!\!\!\! e^{-\cE_2(0)/2} 
 =
\sum_{t=0} \sum_{v=t} \frac{(2d')!}{t!(2d' - v)!} u^t [z^{v-t}]  \frac{\SS(2uz)^t}{\SS(uz)^{2d'+1}}
 \EE_{2t-2d'}(uz).
\end{multline*}
\end{itemize}

\subsection{Proof of Theorem \ref{thm:hz}}
We are now armed to compute the constants $\epsilon_g(d')$.
\begin{proof}
Facts $1-4$ together imply the first equality of the following computation.
\begin{align*}
\epsilon_g(d') &= [u^{2g - 1 + d'}] 
\sum_{t=0} \sum_{v=t} \frac{(2d')!}{t!(2d' - v)!} u^t [z^{v-t}] \SS(uz)^{-2d'-1} \SS(2uz)^t
\cor{
 \EE_{2t-2d'}(uz)}  \qquad \qquad \left(  \cor{\mathcal{E}_a(z)} = 
  \frac{\delta_{a,0}}{\varsigma(z)} \right)\\
&= [u^{2g - 1}] 
 \sum_{w=0} \frac{(2d')!}{d'!(d' - w)!} [z^{w}] \SS(uz)^{-2d'-1} \SS(2uz)^{d'} \frac{(uz)^{-1}}{\SS(uz)}  \qquad \qquad \qquad \qquad \qquad \qquad \left( w+1 = 2g \right) 
 \\
  &= 
 \frac{(2d'-1)!! 2^{d'}}{(d' - 2g + 1)!} [u^{2g}] \left[ \frac{u/2}{\sinh(u/2)}\right]^{2} \left[ \frac{e^u - e^{-u}}{(e^{u/2} - e^{-u/2})^2}  (u/2) \right]^{d'}
 \\
 &= 
 \frac{(2d'-1)!! 2^{d'-2g}}{(d' - 2g + 1)!} [u^{2g}] \left[ \frac{u}{\sinh(u)}\right]^{2} \left[ \frac{u }{\tanh(u)} \right]^{d'}.
\end{align*}
\end{proof}

\appendix

\section{Jucy correspondence and proof of Fact 2}\label{app}
\noindent
The goal of this appendix is to re-prove Fact 2 via Jucys-Murphy correspondence, following \cite{ALS}.
\begin{teo}[Jucys correspondence \cite{J}, 1974] Let $\mu$ be a conjugacy class of $\mathfrak{S}_n$ or equivalently a partition of $n$, and let $\ell(\mu)$ the number of its cycles or parts. Let
$
C_{\mu} = \sum_{g \in \mu} g \in \mathcal{Z}( \Q(\mathfrak{S}_n))
$ 
be the formal sum of all permutations with cycle type $\mu$ in the center of the group algebra. We have
\begin{equation}\label{thm:Jucys}
\sigma_{k}(\J_2, \dots, \J_{n}) = \!\!\! \sum_{\mu : \ell(\mu) = n - k} \!\!\!\!\! C_{\mu} \,\in \mathcal{Z}\Q[\mathfrak{S}_{n}] , \qquad\qquad  k=0, \dots, n-1,
\end{equation}
where $\sigma_k$ is the $k$-th elementary symmetric polynomial and 
$$
\J_k := (1 \; k) + (2 \; k) + \dots + (k\!-\!1 \;\; k) \in \Q[\mathfrak{S}_n]
\qquad 
\qquad
\text{ for }
\;
 k=2,\dots,n,
$$
is the $k$-th Jucys-Murphy element.
By convention $\sigma_0(\vec{\J}) = id \in \mathfrak{S}_n$.
\end{teo}
\noindent
\begin{ese}[Testing Jucys correspondence for $ {\mathcal{Z}\Q[\mathfrak{S}_{4}]}$]

\begin{align*}
\sigma_0(\J_2, \J_3, \J_4) &= id = (1)(2)(3)(4) = \sum_{\mu : \ell(\mu) = 4 - 0 = 4} \!\!\!\!\!\!C_{\mu}
\\
\sigma_1(\J_2, \J_3, \J_4) &= \J_2 + \J_3 + \J_4= (12) + ((13) + (23)) + ((14) + (24) + (34))\\
&=  \text{all transpositions}  = \, (12)(3)(4) + (13)(2)(4) +  \dots = \sum_{\mu : \ell(\mu) = 4 - 1 = 3} \!\!\!\!\!\! C_{\mu}  
\\
\sigma_2(\J_2, \J_3, \J_4) &= \J_2 \J_3 + \J_2\J_4 + \J_3\J_4 =\\
&= (12)\left((13) + (23)\right) + (12)\left((14) + (24) + (34)\right) + \\
&+ \left((13) + (23)\right)\left((14) + (24) + (34)\right)=\\
&= (123)(4) + (132)(4) + (142)(3) + (124)(3) + (12)(34) + (143)(2)+\\
 &+ (14)(23) + (13)(24) + (432)(1) + (134)(2) + (234)(1) = 
\!\!\!\! \sum_{\mu : \ell(\mu) = 4 - 2 = 2} \!\!\!\!\!\! C_{\mu} 
\\
\sigma_3(\J_2, \J_3, \J_4) &= \J_2 \J_3 \J_4 = (12)[(13) + (23)][(14) + (24) + (34)] = \\
&= (1234) + (1243) + (1324) + (1342) + (1423) + (1432) = \!\!\!\! \sum_{\mu : \ell(\mu) = 4 - 3 = 1} \!\!\!\!\!\! C_{\mu} 
\end{align*}
\end{ese}

\subsection{Proof of fact 2} By definition we have
\begin{align*}
h_{g; (2d')}^{\text{Gr.}, 2} &:= \frac{1}{(2d')!}[C_{id}]   C_{(2)^{d'}}\!\! \left( \sum_{\mu : \ell(\mu) = d - (2g - 1 + d')} C_{\mu} \right)\!\!C_{(2d')},
\\
h_{g;  (2d')}^{<, 2} &:= \frac{1}{(2d')!}[C_{id}] C_{(2)^{d'}} \cdot \sigma_{2g - 1 + d'}(\J_2, \dots, \J_{2d'})  C_{(2d')} ,
\end{align*}
where $[C_{id}]$ is the operator that extracts the coefficient of $C_{id}$ from the expression.
Jucys correspondence \ref{thm:Jucys} for $n = 2d'$ and $k = 2g - 1 + d'$ gives
$$
\sigma_{2g - 1 + d'}(\J_2, \dots, \J_{2d'})
 = 
 \!\!\!\!\!\!\!
\sum_{\mu : \ell(\mu) = d - (2g - 1 + d')} \!\!\!\!\!\!C_{\mu},
$$
which immediately implies
$
h_{g;  (2d')}^{<, 2}
=
h_{g; (2d')}^{\text{Gr.}, 2}.
$
The rest of the appendix is dedicated to explain how the definitions above correspond to their respective  geometric meaning.

\subsubsection{One part Grothendieck dessins d'enfant with one $2$-orbifold ramification}
Let $h_{g; (2d')}^{\text{Gr.}, 2}$ be the weighted number of ramified coverings of the Riemann sphere $\mathbb{C}\mathbb{P}^1$ of degree $2d'$ by a genus $g \geq 0$ surface only ramified over $0, 1 , \infty \,\in \,\mathbb{C}\mathbb{P}^1$, whose ramification profiles over $0$ and $\infty$ are given respectively by the cycle types $(2d')$ and $(2, 2, \dots, 2)$, whereas the cycle type $\mu \vdash 2d'$ over $1$ is arbitrary. The Riemann-Hurwitz equation imposes the following restriction on the lenght of $\mu$:
$
 \ell(\mu) =  d' + 1 - 2g = 2d' - (2g - 1 + d').
$
In terms of the group algebra, $h_{g; (2d')}^{\text{Gr.}, 2}$ is defined as
\begin{align*}
h_{g; (2d')}^{\text{Gr.}, 2} := \frac{1}{(2d')!}[C_{id}]   C_{(2)^{d'}} \left( \sum_{\mu : \ell(\mu) = d - (2g - 1 + d')} C_{\mu} \right)C_{(2d')},
\end{align*}
where $[C_{id}]$ is the operator that extracts the coefficient of $C_{id}$ from the expression.
The expression above can be thought as follows. Label the $2d'$ sheets of the covering. Take an unbrached point on the Riemann sphere and compose the three loops around the points $0,1,\infty$. Every loop picks up the monodromy around that point. Since the base curve is a sphere, it is always possible to "pull" the loop on the other side of the sphere and contract it, this means that the product of the monodromies should be equal to the identity. At the level of the group algebra, this corresponds to formally expanding all the $C_{\mu}$ in the expression above, compute all multiplications in $\mathfrak{S}_{2d'}$, and count how many identities one gets this way. The final division by $(2d')!$ accounts for the possible choices of labelling the sheets.
\subsubsection{One part strictly monotone Hurwitz numbers with one $2$-orbifold ramification}
Let $h_{g;  (2d')}^{<, 2}$ be the weighted number of ramified coverings of the Riemann sphere $\mathbb{C}\mathbb{P}^1$ of degree is $2d'$ by a genus $g \geq 0$ surface ramified over $0, \infty \,\in \,\mathbb{C}\mathbb{P}^1$, with ramification profiles given respectively by the cycle types $(2d')$ and $(2, 2, \dots, 2)$, and $m$ further simple ramifications (i.e. these $m$ ramification profiles are given by $m$ transpositions $(x_i \; y_i) \, \in \, \mathfrak{S}_{2d'}$). Let us write them in such a way that $x_i<y_i$, $i=1,\dots,m$. We impose the extra \textit{strictly monotone} condition $y_i < y_{i+1}.$
The Riemann-Hurwitz equation imposes that $m = 2g - 1 + d'$.
In terms of the group algebra, $h_{g;  (2d')}^{<, 2}$ is defined as
\begin{align*}
h_{g;  (2d')}^{<, 2} &:= \frac{1}{(2d')!}[C_{id}] C_{(2)^{d'}} \cdot \sigma_{2g - 1 + d'}(\J_2, \dots, \J_{2d'})  C_{(2d')} .
\end{align*}
Let us explain why the symmetric elementary polynomials $\sigma$ of Jucys-Marphy elements are the right elements to describe these ramified coverings. First of all, every Jucys-Murphy element is a formal sum of transpositions, and $\sigma_{2g - 1 + d'}$ is homogenous of degree $2g - 1 + d'$, therefore every summand in its expansion is a product of exactly $2g - 1 + d'$ transpositions, as required. Let us see how the strict monotonicity condition is guaranteed. The Jucys-Murphy element $\J_k$ is the formal sum of all those ramification profiles whose greatest sheet label is equal to $k$, therefore, if $\J_k$ appears as $i$-th factor in some summand of $\sigma$, this translates to $y_i = k$ in the corresponding covering. At this point, the definition of $\sigma$ grants that the summands produced are all and only the ones satisfying the monotonicity conditions in the corresponding coverings.

\subsection*{Acknowledgements}

This short note was inspired by the talk of N.~Do at the conference \emph{MoSCATR VI: 6th Workshop on Combinatorics of Moduli Spaces, Cluster Algebras, and Topological Recursion} in June 2018 about Hurwitz problems and related one-point recursions. The author is grateful to A.~Alexandrov, R.~Kramer, and S.~Shadrin for useful discussions. The work of the author is supported by the Max-Planck-Gesellschaft.

\end{document}